\documentclass[12pt]{article}
\usepackage{amsmath,amssymb,amsthm,xspace,amscd}
\newtheorem{definition}{Definition}
\newtheorem{remark}{Remark}
\newtheorem{theorem}{Theorem}
\newtheorem{lemma}{Lemma}
\newtheorem{proposition}{Proposition}
\newtheorem{corollary}{Corollary}

\newcommand{\N}{\ensuremath{\mathbb{N}}\xspace}
\newcommand{\Z}{\ensuremath{\mathbb{Z}}\xspace}
\def\bD{{\bf D}}
\def\mod{{\rm mod}}
\def\pro{{\rm pro}}
\def\v{\vector}
\def\bu{\bullet}
\def\a{\alpha}
\def\b{\beta}
\def\d{\delta}
\def\f{\varphi}
\def\g{\gamma}

\begin{document}
\title{Derived categories of Schur algebras}
\author{Viktor Bekkert\thanks{
The first author was supported by FAPESP (Grant 98/14538-0)}
 and  Vyacheslav Futorny\thanks{
\noindent
Regular Associate of the ICTP
\hfil \break
}}

\date{}
\maketitle

\begin{abstract}
In this paper we classify the derived tame Schur 
and infinitesimal Schur algebras 
and describe indecomposable objects in their derived categories.

\end{abstract}

\section{Introduction}\label{s1}

Let $A$ be a finite-dimensional algebra 
of the form $kQ/I$ over an algebraically closed field $k$,  
$A-\mod$ be the category of left $A$-modules and let 
$\bD^b(A)$ be the bounded derived category of the category $A-\mod$. 
The category  $\bD^b(A)$ is well-understood for some classes of algebras $A$.
For example, the description of indecomposable objects of $\bD^{b}(A)$
is well-known for hereditary algebras of finite and tame type \cite{H}
and for the tubular algebras \cite{HR}. 

We say that $A$ is derived tame (see \cite{GK}) 
if for each sequence ${\bf n}=(n_i)_{i\in \Z}$
of positive integers the indecomposable complexes in $\bD^{b}(A)$ 
of cohomology dimension ${\bf n}$ can be parametrized using only one
continuous parameter.

The goal of this paper is to prove the following theorem.

\begin{theorem} Let $p$ be the characteristic of $k$.
\begin{itemize}
\item[(i)] Let $S=S(n,d)$  be a Schur algebra. 
Then $S$ is derived tame if and only if one of the following holds:

a) $p>d$;

b) $p=2$, $n=2$, $d=3$;

c) $p=2$, $n=2$, $d=5$ or $d=7$;

d) $n=2$, $p\leq d< 2p$ and $d\neq 3$ if $p=2$; 

e) $p=2$,  $n\geq 3$, $d=2$ or $d=3$;

f) $p=3$, $n=3$, $d=4$ or $d=5$.

\item[(ii)] Let $S=S(n,d)_r$ be an infinitesimal Schur algebra. Then $S$ is derived tame 
if and only if one of the following holds: 

a) $p > 2$, $n\geq 2$, $d<p$;

b) $p=2$, $n=2$, $d=3$;

c) $p=2$, $n=2$, $d$ odd, $r=1$;

d) $p=2$, $n=2$, $r=2$, $d=5$ or $d=7$;

e) $p=2$, $n=2$, $r=1$, $d=2$;

f) $n=2$, $r=1$, $2<p\leq d<2p$.

\end{itemize}

\end{theorem}

The structure of the paper is the following. 
In Section 2 the definition of Schur algebras 
and some preliminary results from  \cite{E} and \cite{DEMN}  are given.
In Section 3 we  show that the problem of classification of all
indecomposable objects in $\bD^{b}(A)$ can be reduced to the problem of classification 
of indecomposable objects in the category $\frak p(A)$, some subcategory
of the category of bounded projective complexes ${\bf C}^{b}(A-\pro)$
and we prove  Theorem 1. In Section 4 
the description of all indecomposable objects in $\bD^{b}(A)$
for derived tame Schur and infinitesimal Schur algebras is given.

\section{Preliminaries}

Let $V$ be an $n$-dimensional vector space over an algebraically closed field of characteristic 
$p$. Then the symmetric group $S_d$ has a natural action on $V^{\otimes d}$ which makes it a module 
for the group algebra of $S_d$. The endomorphism ring of the 
module $V^{\otimes d}$ is the Schur algebra 
$S(n, d)$. There is an equivalence between the category of modules for $S(n, d)$ and the category 
of polynomial representations of $GL(n)$ which are homogeneous of degree $d$ \cite{G}. Schur 
algebras are quasi-hereditary \cite{P} like the algebras corresponding to the blocks of category 
${\cal O}$ \cite{BGG}. All Schur algebras of finite representation type, 
i.e., those algebras that have only finitely many non-isomorphic
 indecomposable modules, were classified in \cite{E}. 
 
The Schur algebra is of finite representation type if one of the following holds:
\medskip

 a) $n=2$ and $d<p^2$ or $n\geq 3$ and $d<2p$;

 b) $p=n=2$ and $d=5$ or $d=7$. 
\medskip

Moreover, the corresponding  blocks
  are Morita equivalent 
to one of the algebras $A_m$ below:
\bigskip

\begin{picture}(400.00,50.00)
\put(0.00,16.00){$A_m:$}
\put(0.00,2.00){$_{m\geq 1}$}
\put(260.00,25.00){$\a_i\a_{i+1}=0, \b_{i+1}\b_i=0,$}
\put(260.00,5.00){$\b_i\a_i=\a_{i+1}\b_{i+1}, \a_1\b_1=0.$}
\put(30.00,16.00){$\bu$}
\put(30.00,26.00){$1$}
\put(80.00,16.00){$\bu$}
\put(80.00,26.00){$2$}
\put(130.00,16.00){$\bu$}
\put(130.00,26.00){$3$}
\put(180.00,16.00){$\bu$}
\put(166.00,26.00){${m-1}$}
\put(230.00,16.00){$\bu$}
\put(230.00,26.00){$m$}
\put(150.00,16.00){$\cdots$}
\put(52.00,29.00){$\a_1$}
\put(40.00,23.00){\v(1,0){35.00}}
\put(52.00,0.00){$\b_1$}
\put(75.00,14.00){\v(-1,0){35.00}}
\put(102.00,29.00){$\a_2$}
\put(90.00,23.00){\v(1,0){35.00}}
\put(102.00,0.00){$\b_2$}
\put(125.00,14.00){\v(-1,0){35.00}}
\put(202.00,29.00){$\a_{m-1}$}
\put(190.00,23.00){\v(1,0){35.00}}
\put(202.00,0.00){$\b_{m-1}$}
\put(225.00,14.00){\v(-1,0){35.00}}
\end{picture}

\bigskip

According  to \cite{D}, infinite representation type Schur 
algebras  either have tame or wild 
representation type. The algebra is of tame 
representation type if it has at most finitely many one-parameter
families of indecomposable modules in each dimension. Otherwise, 
 the algebra is wild.  The Schur 
algebras of tame representation type were classified 
in \cite{DEMN} and are covered by  the following
cases:
\medskip
 a) $p=n=3$ and $d=7$ or $d=8$;

 b) $p=3$, $n=2$, $d=9$ or $d=10$ or $d=11$;

 c) $p=n=2$, $d=4$ or $d=9$.
\medskip

The corresponding non-semisimple blocks are Morita equivalent to the 
 algebras given by the following quivers 
with relations.  

\bigskip
\noindent 
1. For algebras $S(2,4)$, $S(2,9)$, $p=2$

\begin{picture}(400.00,50.00)
\put(0.00,16.00){$D_3:$}
\put(160.00,25.00){$\a_1\b_1=\a_2\b_2=0,$}
\put(160.00,5.00){$\a_1\b_2\a_2=\b_2\a_2\b_1=0.$}
\put(30.00,16.00){$\bu$}
\put(30.00,26.00){$1$}
\put(80.00,16.00){$\bu$}
\put(80.00,26.00){$2$}
\put(130.00,16.00){$\bu$}
\put(130.00,26.00){$3$}
\put(52.00,29.00){$\a_1$}
\put(40.00,23.00){\v(1,0){35.00}}
\put(52.00,0.00){$\b_1$}
\put(75.00,14.00){\v(-1,0){35.00}}
\put(102.00,29.00){$\b_2$}
\put(90.00,23.00){\v(1,0){35.00}}
\put(102.00,0.00){$\a_2$}
\put(125.00,14.00){\v(-1,0){35.00}}
\end{picture}

\bigskip
\noindent
2. For algebras $S(2,9)$, $S(2,10)$, $S(2,11)$, $p=3$

\begin{picture}(114.00,100.00)
\put(0.00,14.00){$D_4:$}
\put(160.00,44.00){$\a_1\b_1=\a_2\b_2=\a_3\b_1=\a_3\b_2=0,$}
\put(160.00,24.00){$\a_1\b_3=\a_2\b_3=0, \b_2\a_2=\b_3\a_3,$ }
\put(160.00,4.00){$\a_1\b_2\a_2=\b_2\a_2\b_1=0.$ }
\put(30.00,14.00){$\bu$}
\put(30.00,24.00){$1$}
\put(80.00,14.00){$\bu$}
\put(80.00,00.00){$0$}
\put(80.00,64.00){$\bu$}
\put(80.00,74.00){$2$}
\put(130.00,14.00){$\bu$}
\put(130.00,24.00){$3$}
\put(50.00,27.00){$\a_1$}
\put(50.00,0.00){$\b_1$}
\put(63.00,42.00){$\b_2$}
\put(94.00,42.00){$\a_2$}
\put(100.00,27.00){$\b_3$}
\put(100.00,0.00){$\a_3$}
\put(40.00,21.00){\v(1,0){35.00}}
\put(75.00,12.00){\v(-1,0){35.00}}
\put(90.00,21.00){\v(1,0){35.00}}
\put(125.00,12.00){\v(-1,0){35.00}}
\put(88.00,60.00){\v(0,-1){35.00}}
\put(78.00,25.00){\v(0,1){35.00}}
\end{picture}

\bigskip
\noindent 
3. For the algebra $S(3,7)$,  $p=3$

\begin{picture}(400.00,70.00)
\put(0.00,16.00){$R_4:$}
\put(210.00,25.00){$\a_1\b_1=\a_1\a_2=\b_2\b_1=0,$}
\put(210.00,5.00){$\b_1\a_1=\a_2\b_2, \b_2\a_2=\a_3\b_3.$}
\put(30.00,16.00){$\bu$}
\put(30.00,26.00){$1$}
\put(80.00,16.00){$\bu$}
\put(80.00,26.00){$2$}
\put(130.00,16.00){$\bu$}
\put(130.00,26.00){$3$}
\put(180.00,16.00){$\bu$}
\put(180.00,26.00){$4$}
\put(52.00,29.00){$\a_1$}
\put(40.00,23.00){\v(1,0){35.00}}
\put(52.00,0.00){$\b_1$}
\put(75.00,14.00){\v(-1,0){35.00}}
\put(102.00,29.00){$\a_2$}
\put(90.00,23.00){\v(1,0){35.00}}
\put(102.00,0.00){$\b_2$}
\put(125.00,14.00){\v(-1,0){35.00}}
\put(152.00,29.00){$\a_3$}
\put(140.00,23.00){\v(1,0){35.00}}
\put(152.00,0.00){$\b_3$}
\put(175.00,14.00){\v(-1,0){35.00}}
\end{picture}

\bigskip
\noindent 
4. For the algebra $S(3,8)$, $p=3$

\begin{picture}(114.00,85.00)
\put(0.00,14.00){$H_4:$}
\put(160.00,44.00){$\a_1\b_1=\a_1\b_2=\a_1\b_3=0,$}
\put(160.00,24.00){$\a_3\b_1=\a_3\b_3=\a_2\b_1=0,$ }
\put(160.00,4.00){$\b_2\a_2=\b_1\a_1 + \b_3\a_3.$ }
\put(30.00,14.00){$\bu$}
\put(30.00,24.00){$1$}
\put(80.00,14.00){$\bu$}
\put(80.00,00.00){$0$}
\put(80.00,64.00){$\bu$}
\put(80.00,74.00){$2$}
\put(130.00,14.00){$\bu$}
\put(130.00,24.00){$3$}
\put(50.00,27.00){$\a_1$}
\put(50.00,0.00){$\b_1$}
\put(63.00,42.00){$\b_2$}
\put(94.00,42.00){$\a_2$}
\put(100.00,27.00){$\b_3$}
\put(100.00,0.00){$\a_3$}
\put(40.00,21.00){\v(1,0){35.00}}
\put(75.00,12.00){\v(-1,0){35.00}}
\put(90.00,21.00){\v(1,0){35.00}}
\put(125.00,12.00){\v(-1,0){35.00}}
\put(88.00,60.00){\v(0,-1){35.00}}
\put(78.00,25.00){\v(0,1){35.00}}
\end{picture}

\bigskip

In \cite{DNP1} the authors introduced certain subalgebras $S(n,d)_r$ which they called 
{\em infinitesimal Schur algebras}. The representation theory of these algebras is connected 
to  the theory of polynomial representations of the group scheme $G_rT$. Here $G=GL(n)$ over 
an algebraically closed field of characteristic $p>0$, $G_rT$ is the inverse 
 image of the diagonal torus $T\subset G$ under the $r$th iteration of the Frobenius 
 morphism. The infinitesimal Schur algebras $S(n,d)_r$ of finite 
 representation type were classified in \cite{DNP2}. They belong to
 the following cases: 

\medskip
a) $p\geq 2$, $n\geq 3$, $d<2p$, $r\geq 2$;

b) $p\geq 2$, $n\geq 3$, $d < p$, $r=1$;

c) $p=3$, $n=3$, $r=1$, $d=4$ or $d=5$;

d) $p=2$, $n=3$, $r=1$, $d=2$ or $d=3$;

e) $p\geq 2$, $n=2$, $d<p^2$, $r\geq 2$;

f) $p=2$, $n=2$, $r\geq 3$, $d=5$ or $d=7$;

g) $p=2$, $n=2$, $r=2$, $d$ odd;

h) $p\geq 2$, $n=2$, $r=1$.

\medskip
The infinitesimal Schur algebras in cases a), b), e), f) 
 coincide with the corresponding Schur algebras. In all other cases  
 the corresponding non-semisimple blocks are Morita 
 equivalent to the algebras given by the following quivers with relations.

\bigskip
\noindent 1. For algebras $S(3,4)_1$, $S(3,5)_1$, $p=3$ and  
$S(3,2)_1$, $S(3,3)_1$, $p=2$

\begin{picture}(114.00,100.00)
\put(0.00,24.00){$G:$}
\put(160.00,34.00){$\b_1\a_1=\b_2\a_2=\b_3\a_3,$}
\put(160.00,14.00){all other products $=0$. }
\put(30.00,24.00){$\bu$}
\put(30.00,34.00){$1$}
\put(80.00,24.00){$\bu$}
\put(80.00,8.00){$0$}
\put(80.00,74.00){$\bu$}
\put(80.00,84.00){$2$}
\put(130.00,24.00){$\bu$}
\put(130.00,34.00){$3$}
\put(50.00,37.00){$\a_1$}
\put(50.00,10.00){$\b_1$}
\put(63.00,52.00){$\b_2$}
\put(94.00,52.00){$\a_2$}
\put(100.00,37.00){$\b_3$}
\put(100.00,10.00){$\a_3$}
\put(40.00,31.00){\v(1,0){35.00}}
\put(75.00,22.00){\v(-1,0){35.00}}
\put(90.00,31.00){\v(1,0){35.00}}
\put(125.00,22.00){\v(-1,0){35.00}}
\put(88.00,70.00){\v(0,-1){35.00}}
\put(78.00,35.00){\v(0,1){35.00}}
\end{picture}

\bigskip
\noindent 
2. For algebras $S(2,d)_2$, $d$ odd, $p=2$ and $S(2,d)_1$, $p\geq 2$

\begin{picture}(400.00,80.00)
\put(0.00,26.00){$F_m:$}
\put(0.00,10.00){$_{m=2n+1}$}
\put(0.00,2.00){$_{n>0}$}
\put(260.00,35.00){$\a_i\a_{i+1}=0, \b_{i+1}\b_i=0, \a_{1}\b_{1}=0,$}
\put(260.00,15.00){$\b_i\a_i=\a_{i+1}\b_{i+1}, \b_{m-1}\a_{m-1}=0.$}
\put(30.00,26.00){$\bu$}
\put(30.00,36.00){$1$}
\put(80.00,26.00){$\bu$}
\put(80.00,36.00){$2$}
\put(130.00,26.00){$\bu$}
\put(130.00,36.00){$3$}
\put(180.00,26.00){$\bu$}
\put(166.00,36.00){$m-1$}
\put(230.00,26.00){$\bu$}
\put(230.00,36.00){$m$}
\put(150.00,26.00){$\cdots$}
\put(52.00,39.00){$\a_1$}
\put(40.00,33.00){\v(1,0){35.00}}
\put(52.00,10.00){$\b_1$}
\put(75.00,24.00){\v(-1,0){35.00}}
\put(102.00,39.00){$\a_2$}
\put(90.00,33.00){\v(1,0){35.00}}
\put(102.00,10.00){$\b_2$}
\put(125.00,24.00){\v(-1,0){35.00}}
\put(202.00,39.00){$\a_{m-1}$}
\put(190.00,33.00){\v(1,0){35.00}}
\put(202.00,10.00){$\b_{m-1}$}
\put(225.00,24.00){\v(-1,0){35.00}}
\end{picture}

\bigskip

The infinitesimal Schur algebras of tame representation type were classified in 
 \cite{DEMN}. They belong to  the following classes:

\medskip
a) $p\geq 5$, $n=3$, $p\leq d\leq 2p-1$, $r=1$;

b) $p=3$, $n=3$, $d=3$, $r=1$;

c) $p=3$, $n=4$, $r=1$, $d=3$ or $d=4$ or $d=5$;

d) $p=3$, $n=3$, $r\geq 2$, $d=7$ or $d=8$;

e) $p=3$, $n=2$, $r\geq 3$, $d=9$ or $d=10$ or $d=11$;

f) $p=2$, $n=4$, $r=1$, $d=2$ or $d=3$;

g) $p=2$, $n=2$, $d=4$, $r\geq 2$;

h) $p=2$, $n=2$, $d=9$, $r\geq 3$.

\medskip
The algebras in e) have non-semisimple blocks of type $D_4$. The algebra 
$S(3,7)_r$, $r\geq 2$, $p=3$ has a block of type $R_4$. The algebra $S(3,8)_r$, 
$r\geq 2$, $p=3$ has a block of type $H_4$. All algebras in g) with $r>2$ and 
all algebras in h) with $r>3$ have blocks of type $D_3$. Other algebras have 
the blocks which are Morita equivalent to the following quivers with relations:

\bigskip
\noindent 
1. For algebras $S(4,3)_1$, $S(4,4)_1$, $S(4,5)_1$, $p=3$ and 
$S(4,2)_1$, $S(4,3)_1$, $p=2$
  
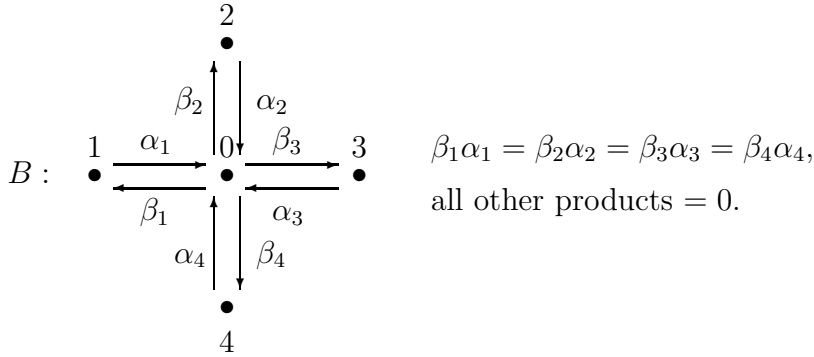
\begin{picture}(114.00,150.00)
\put(0.00,70.00){$B:$}
\put(160.00,80.00){$\b_1\a_1=\b_2\a_2=\b_3\a_3=\b_4\a_4,$}
\put(160.00,60.00){all other products $=0$.}
\put(30.00,70.00){$\bu$}
\put(30.00,80.00){$1$}
\put(80.00,70.00){$\bu$}
\put(80.00,80.00){$0$}
\put(80.00,20.00){$\bu$}
\put(80.00,6.00){$4$}
\put(80.00,120.00){$\bu$}
\put(80.00,130.00){$2$}
\put(130.00,70.00){$\bu$}
\put(130.00,80.00){$3$}
\put(50.00,83.00){$\a_1$}
\put(50.00,56.00){$\b_1$}
\put(63.00,98.00){$\b_2$}
\put(94.00,98.00){$\a_2$}
\put(100.00,83.00){$\b_3$}
\put(100.00,56.00){$\a_3$}
\put(94.00,40.00){$\b_4$}
\put(63.00,40.00){$\a_4$}
\put(40.00,77.00){\v(1,0){35.00}}
\put(75.00,68.00){\v(-1,0){35.00}}
\put(90.00,77.00){\v(1,0){35.00}}
\put(125.00,68.00){\v(-1,0){35.00}}
\put(88.00,65.00){\v(0,-1){35.00}}
\put(78.00,30.00){\v(0,1){35.00}}
\put(88.00,116.00){\v(0,-1){35.00}}
\put(78.00,81.00){\v(0,1){35.00}}
\end{picture}

\bigskip
\noindent 
2. For algebras $S(3, d)_1$, $p\geq 5$, $p\leq d\leq 2p-1$ and 
$S(3,3)_1$, $p=3$

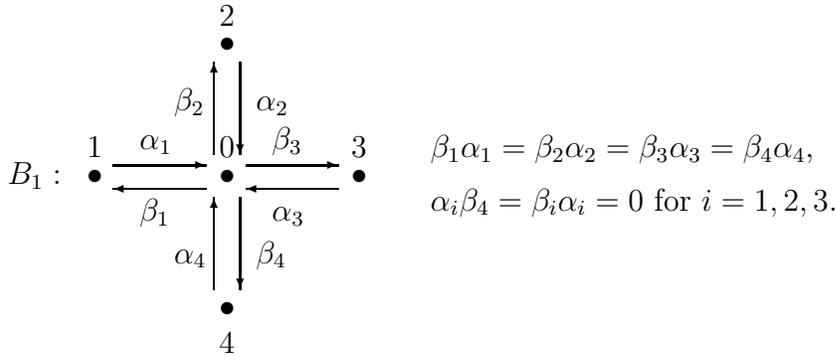
\begin{picture}(114.00,150.00)
\put(0.00,70.00){$B_1:$}
\put(160.00,80.00){$\b_1\a_1=\b_2\a_2=\b_3\a_3=\b_4\a_4,$}
\put(160.00,60.00){$\a_i\b_4=\b_i\a_i=0$ for $i=1,2,3$. }
\put(30.00,70.00){$\bu$}
\put(30.00,80.00){$1$}
\put(80.00,70.00){$\bu$}
\put(80.00,80.00){$0$}
\put(80.00,20.00){$\bu$}
\put(80.00,6.00){$4$}
\put(80.00,120.00){$\bu$}
\put(80.00,130.00){$2$}
\put(130.00,70.00){$\bu$}
\put(130.00,80.00){$3$}
\put(50.00,83.00){$\a_1$}
\put(50.00,56.00){$\b_1$}
\put(63.00,98.00){$\b_2$}
\put(94.00,98.00){$\a_2$}
\put(100.00,83.00){$\b_3$}
\put(100.00,56.00){$\a_3$}
\put(94.00,40.00){$\b_4$}
\put(63.00,40.00){$\a_4$}
\put(40.00,77.00){\v(1,0){35.00}}
\put(75.00,68.00){\v(-1,0){35.00}}
\put(90.00,77.00){\v(1,0){35.00}}
\put(125.00,68.00){\v(-1,0){35.00}}
\put(88.00,65.00){\v(0,-1){35.00}}
\put(78.00,30.00){\v(0,1){35.00}}
\put(88.00,116.00){\v(0,-1){35.00}}
\put(78.00,81.00){\v(0,1){35.00}}
\end{picture}

\bigskip
\noindent 3. For algebras $S(2,4)_2$,  $S(2,9)_3$, $p=2$

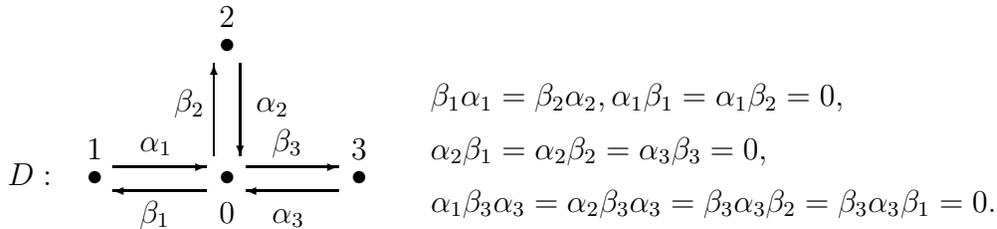
\begin{picture}(114.00,100.00)
\put(0.00,14.00){$D:$}
\put(160.00,44.00){$\b_1\a_1=\b_2\a_2, \a_1\b_1=\a_1\b_2=0,$}
\put(160.00,24.00){$\a_2\b_1=\a_2\b_2=\a_3\b_3=0,$ }
\put(160.00,4.00){$\a_1\b_3\a_3=\a_2\b_3\a_3=\b_3\a_3\b_2=\b_3\a_3\b_1=0.$ }
\put(30.00,14.00){$\bu$}
\put(30.00,24.00){$1$}
\put(80.00,14.00){$\bu$}
\put(80.00,00.00){$0$}
\put(80.00,64.00){$\bu$}
\put(80.00,74.00){$2$}
\put(130.00,14.00){$\bu$}
\put(130.00,24.00){$3$}
\put(50.00,27.00){$\a_1$}
\put(50.00,0.00){$\b_1$}
\put(63.00,42.00){$\b_2$}
\put(94.00,42.00){$\a_2$}
\put(100.00,27.00){$\b_3$}
\put(100.00,0.00){$\a_3$}
\put(40.00,21.00){\v(1,0){35.00}}
\put(75.00,12.00){\v(-1,0){35.00}}
\put(90.00,21.00){\v(1,0){35.00}}
\put(125.00,12.00){\v(-1,0){35.00}}
\put(88.00,60.00){\v(0,-1){35.00}}
\put(78.00,25.00){\v(0,1){35.00}}
\end{picture}

\bigskip

\section{Derived representation type}

Let $A$ be a finite-dimensional algebra of the form $kQ/I$ 
over 
an algebraically closed  field $k$,  
$A-\mod$ be the category of left $A$-modules. 
We will follow in general the notation and terminology of \cite{Ri} and 
\cite{H}.

Given $A$,  we denote by  ${\bf D}(A)$ (resp., 
${\bf D}^{-}(A)$ or
${\bf D}^{b}(A)$) the derived category of $A-\mod$
(resp., 
the derived category of right bounded complexes of $A-\mod$
or
the derived category of bounded complexes of $A-\mod$); by
${\bf C}^{b}(A-\pro)$ (resp., ${\bf C}^{-}(A-\pro)$ or
${\bf C}^{-,b}(A-\pro)$) the category of bounded
projective complexes (resp., 
of right bounded projective complexes or
of right bounded projective complexes with
bounded cohomology (that is, complexes of projective modules with the
property
that the cohomology groups are non zero only at a finite number of places));
and by ${\bf K}^{b}(A-\pro)$ 
(resp., ${\bf K}^{-}(A-\pro)$ or 
${\bf K}^{-,b}(A-\pro)$) the
corresponding
homotopy categories.

We identify the homotopy category ${\bf K}^{b}(A-\pro)$ with the full subcategory
of perfect complexes in ${\bf D}^{b}(A)$.
Recall that a complex is perfect if it is isomorphic to a bounded
complex of finitely generated projective $A$-modules.

We will also use the following notations.
By ${\frak p} (A)$ we denote the full
subcategory
of ${\bf C}^{b}(A-\pro)$ 
defined by the projective
complexes such that the image of every differential map is contained
in the radical of the corresponding projective module.
Since any projective complex
is the sum of one
complex with this property and two complex where,
alternativelly, all differential
maps are 0's or isomorphisms (which is, hence, isomorphic to the zero
object in the derived category) we can always assume that we reduce
ourselves
to consider projective complexes of this form.
 
It is well known that $\bD^{b}(A)$   is
equivalent to ${\bf K}^{-,b}(A-\pro)$ \, 
(see, for example, \cite{KZ}, Prop. 6.3.1 and \cite{Har}). 

\begin{proposition}\cite{Har} ${\bf D}^{-}(A)$ is equivalent 
to ${\bf K}^{-}(A-\pro)$.
The image
of \, ${\bf D}^{b}(A)$  under this equivalence is \, ${\bf K}^{-,b}(A-\pro)$. 
\end{proposition}

Given  $M^{\bullet} \in {\bf D}^{b}(A)$  we denote by
$ P^{\bullet }_{M^{\bullet}} $  the
projective resolution of $M^{\bullet}$ (see \cite{KZ})
and by $H^{i}(M^{\bullet})$ the $i$-th cohomology module.
 
We call a category ${\cal C}$ {\em basic} if it satisfies the
following conditions:
\begin{itemize}
\item all its objects are pairwise non-isomorphic;
\item for each object $x$ there are no non-trivial
idempotents in ${\cal C}(x,x)$.
\end{itemize}

A full subcategory ${\cal S}\subset {\cal C}$ is
called a {\em skeleton} of ${\cal C}$ if it is basic and each object
$x\in {\cal C}$ is isomorphic to a direct summand of a
(finite) direct sum of some objects of ${\cal S}$.
It is evident that if  ${\cal C}$ is a category with
unique direct decomposition property, then it has a
skeleton and the last one is unique up to isomorphism.
We will denote it by ${\bf Sk}\,{\cal C}$ and the set of its
objects by ${\bf Ver}\,{\cal C}$.

In order to simplify our exposition, let us introduce two
easy constructions, as follows.

For $P^{\bullet}\in {\bf C}^{-,b}(A-\pro)\setminus {\bf C}^{b}(A-\pro)$, let $s$ be
the maximal  number such that $P^{s} \not = 0$ and 
 $H^{i}(P^{\bullet}) = 0$ for $i \leq s$.
Then,
$\alpha(P^{\bullet})^{\bullet}$ denotes 
the {\em brutal truncation} of $P^{\bullet}$ below $s$
(see \cite{W}),
i.e. the complex given by
$$
\alpha(P^{\bullet})^{i} =\left\{ \begin{array}{ll}
P^{i} & \mbox{, if $i\geq s$;}\\
0 & \mbox{, otherwise,}
\end{array}
\right.
$$

$$
\partial^{i}_{\alpha(P^{\bullet})^{\bullet}} =\left\{ \begin{array}{ll}
\partial^{i}_{P^{\bullet}} & \mbox{, if $i\geq s$;}\\
0 & \mbox{, otherwise.}
\end{array}
\right.
$$
\medskip \noindent For $P^{\bullet}\not=0^{\bullet} \in {\bf C}^{b}(A-\pro)$,
let $t$ be
the maximal  number such that $P^{i} = 0$ for $i < t$. Then,
$\beta(P^{\bullet})^{\bullet}$ denotes 
the {\em (good) truncation} of $P^{\bullet}$ below $t$
(see \cite{W}),
i.e. 
the complex given by
$$
\beta(P^{\bullet})^{i} =\left\{ \begin{array}{lll}
P^{i} & \mbox{, if $i\geq t$;}\\
{\rm Ker}\, \partial_{P^{\bullet}}^{t} & \mbox{, if $i=t-1$;}\\
0 & \mbox{, otherwise,}
\end{array}
\right.
$$

$$
\partial^{i}_{\beta(P^{\bullet})^{\bullet}} =\left\{ \begin{array}{lll}
\partial^{i}_{P^{\bullet}} & \mbox{, if $i\geq t$;}\\
i_{{\rm Ker}\, \partial_{P^{\bullet}}^{t}} & \mbox{, if $i=t-1$;}\\
0 & \mbox{, otherwise,}
\end{array}
\right.
$$

\noindent where $i_{{\rm Ker}\, \partial^{t}_{P^{\bullet}}} $ is the obvious
inclusion.

\begin{lemma} Let $M^{\bullet} 
\in {\bf K}^{-,b}(A-\pro)\setminus {\bf K}^{b}(A-\pro)$ be an indecomposable. Then 
$N^{\bullet}=\beta(\alpha(M^{\bullet})^{\bullet})^{\bullet}$ is also indecomposable 
in $\bD^{b}(A) \,$ and \[ M^{\bullet} \cong 
P^{\bullet}_{N^{\bullet}} . \] 
\end{lemma}

\begin{proof}
Obvious.
\end{proof}

\begin{lemma}
There exist skeletons 
${\bf Sk}\, \frak p(A)$ and
${\bf Sk} \, {\bf K}^{b}(A-\pro)$ of  $\frak p(A)$
and ${\bf K}^{b}(A-\pro)$, respectively, such that
${\bf Ver}\,{\frak p}(A)=$ 
${\bf Ver}\,{\bf K}^{b}(A-\pro)$.
\end{lemma}

\begin{proof}
Obvious.   
\end{proof} 

Let $\overline{{\cal X}(A)}=$ $\{\,M^{\bullet}\in {\bf Ver}\,{\frak p}(A)\,|$
$P^{\bullet}_{\b(M^{\bullet})^{\bullet}}\not\in {\bf K}^{b}(A-\pro)\,\}$.
Let $\cong_{\cal X}$ be the equivalence relation on the set 
$\overline{{\cal X}(A)}$ defined by $M^{\bullet}\cong_{\cal X} N^{\bullet}$
if and only if $P^{\bullet}_{\b(M^{\bullet})^{\bullet}}\cong$
$P^{\bullet}_{\b(N^{\bullet})^{\bullet}}$ in ${\bf K}^{-,b}(A-\pro)$.
We use the notation ${\cal X}(A)$ for a fixed set of representatives of the quotient
set $\overline{{\cal X}(A)}$ over the equivalence relation $\cong_{\cal X}$.

>From Lemmas 1 and 2 we obtain the following

\begin{corollary} There exist skeletons ${\bf Sk}\, {\bf D}^{b}(A)$ and
${\bf Sk}\, {\frak p}(A)$ of ${\bf D}^{b}(A)$ and ${\frak p}(A)$, respectively, such
that ${\bf Ver}\, {\bf D}^{b}(A) = {\bf Ver}\,{\frak p}(A) 
\cup \{\beta(M^{\bullet})^{\bullet}\,|\, M^{\bullet} \in {\cal X}(A)\}.$ 
\end{corollary}

\begin{remark}
If $A$ has a finite global dimension then ${\cal X}(A)=\emptyset$
and ${\bf Ver}\, {\bf D}^{b}(A) = {\bf Ver}\,{\frak p}(A)$.
\end{remark}

Let $T$ be the translation functor 
${\bf D}(A)\rightarrow {\bf D}(A)$. 
By analogy 
to \cite{D} we will use the following definitions.

\begin{definition}
Let $k$ be an algebraically closed field and
$A$ be a finite-dimensional $k$-algebra. Then
\begin{itemize}
\item $A$ is called {\em derived wild} if there exists
a complex of  $A-k\left< x,y\right>$-
bimodules $M^{\bullet}$ such that each $M^{i}$ is free and
of finite rank as right $k\left< x,y\right>$-module and such that the
functor $M\otimes_{k\left< x,y\right> }-$ preserves indecomposability
and isomorphism classes.
\item $A$ is called {\em derived tame}  ({\em see \cite{GK}}) if,
for each cohomology dimension vector
$(d_{i})_{i \in \Z }$,
there exist a localization $R=k[x]_f$ with respect to some
$f\in k[x]$ and  a finite number of  bounded complexes of $A-R$-bimodules
$C_{1}^{\bullet},\cdots , C_{n}^{\bullet}$ such that each $C_{j}^{i}$
is free and of finite rank as  right $R$-module and such
that every indecomposable $X^{\bullet} \in {\bf D}^{b}(A)$  with
${\rm dim} \, H^{i}(X^{\bullet} )=d_{i}$
is isomorphic to $C_{j}^{\bullet} \otimes_{R}  S$ for some
$j$ and
some simple $R$-module $S$.
\item $A$ is called {\em derived  discrete (see \cite{V})} if
for every cohomology
dimension vector
$(d_{i})_{i \in \Z}$,
we  have up to isomorphism a finite number of indecomposables
$X^{\bullet} \in {\bf D}^{b}(A)$
with ${\rm dim} \, H^{i}(X^{\bullet} )=d_{i}$.
\item
$A$  is called {\em derived finite} if we have a finite number of
indecomposables  $$X_{1}^{\bullet} ,\cdots , X_{n}^{\bullet} \in
{\bf D}^{b}(A)$$
such that every   indecomposable object $X^{\bullet} \in$ $ {\bf D}^{b}(A)$
is isomorphic to  $T^{i}(X_j^{\bullet})$ for some $i\in \Z$
and some $j$.
\end{itemize}
\end{definition}

We denote by ${\bf P}_i$ the indecomposable projective corresponding to
the vertex $i\in Q_0$ and by $p(w)$ the morphism between two indecomposable
projectives corresponding to the path $w$ of $Q$.

\subsection{}

We list below the wild algebras which are used in the proof of the Theorem 1.

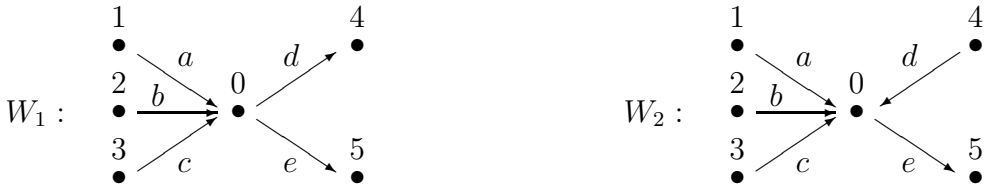
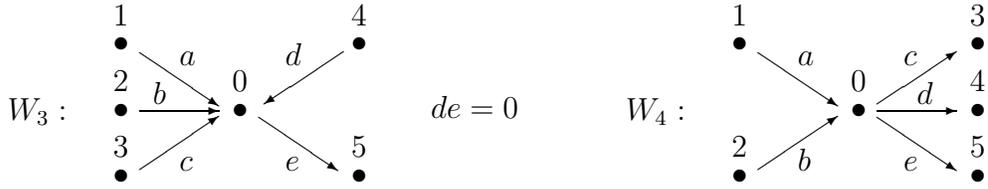
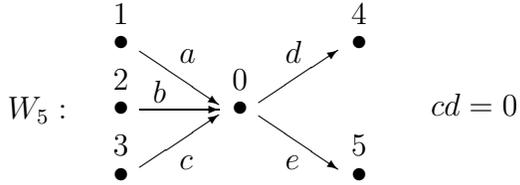
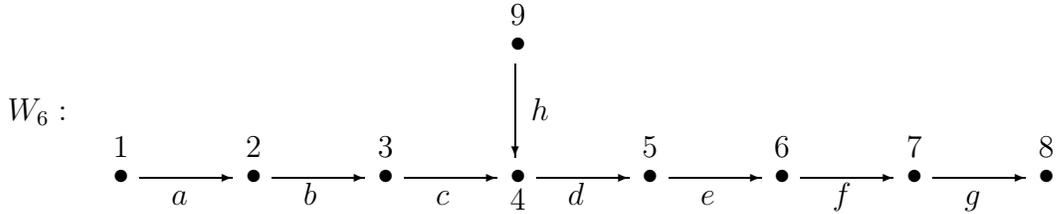
\begin{picture}(230.00,90.00)
\put(0.00,24.00){$W_1:$}
\put(85.00,25.00){$\bu$}
\put(85.00,35.00){$0$}
\put(65.00,45.00){$a$}
\put(65.00,5.00){$c$}
\put(55.00,30.00){$b$}
\put(105.00,45.00){$d$}
\put(105.00,5.00){$e$}
\put(130.00,00.00){$\bu$}
\put(130.00,10.00){$5$}
\put(130.00,50.00){$\bu$}
\put(130.00,60.00){$4$}
\put(95.00,30.00){\v(3,2){30.00}}
\put(95.00,25.00){\v(3,-2){30.00}}
\put(40.00,25.00){$\bu$}
\put(40.00,35.00){$2$}
\put(40.00,0.00){$\bu$}
\put(40.00,10.00){$3$}
\put(40.00,50.00){$\bu$}
\put(40.00,60.00){$1$}
\put(50.00,27.50){\v(1,0){30.00}}
\put(50.00,5.50){\v(3,2){30.00}}
\put(50.00,49.50){\v(3,-2){30.00}}
\end{picture}
\begin{picture}(114.00,90.00)
\put(0.00,24.00){$W_2:$}
\put(85.00,25.00){$\bu$}
\put(85.00,35.00){$0$}
\put(65.00,45.00){$a$}
\put(65.00,5.00){$c$}
\put(55.00,30.00){$b$}
\put(105.00,45.00){$d$}
\put(105.00,5.00){$e$}
\put(130.00,00.00){$\bu$}
\put(130.00,10.00){$5$}
\put(130.00,50.00){$\bu$}
\put(130.00,60.00){$4$}
\put(127.00,50.00){\v(-3,-2){30.00}}
\put(95.00,25.00){\v(3,-2){30.00}}
\put(40.00,25.00){$\bu$}
\put(40.00,35.00){$2$}
\put(40.00,0.00){$\bu$}
\put(40.00,10.00){$3$}
\put(40.00,50.00){$\bu$}
\put(40.00,60.00){$1$}
\put(50.00,27.50){\v(1,0){30.00}}
\put(50.00,5.50){\v(3,2){30.00}}
\put(50.00,49.50){\v(3,-2){30.00}}
\end{picture}

\begin{picture}(230.00,90.00)
\put(0.00,24.00){$W_3:$}
\put(85.00,25.00){$\bu$}
\put(85.00,35.00){$0$}
\put(65.00,45.00){$a$}
\put(65.00,5.00){$c$}
\put(55.00,30.00){$b$}
\put(105.00,45.00){$d$}
\put(105.00,5.00){$e$}
\put(130.00,00.00){$\bu$}
\put(130.00,10.00){$5$}
\put(130.00,50.00){$\bu$}
\put(130.00,60.00){$4$}
\put(127.00,50.00){\v(-3,-2){30.00}}
\put(95.00,25.00){\v(3,-2){30.00}}
\put(40.00,25.00){$\bu$}
\put(40.00,35.00){$2$}
\put(40.00,0.00){$\bu$}
\put(40.00,10.00){$3$}
\put(160.00,25.00){$de=0$}
\put(40.00,50.00){$\bu$}
\put(40.00,60.00){$1$}
\put(50.00,27.50){\v(1,0){30.00}}
\put(50.00,5.50){\v(3,2){30.00}}
\put(50.00,49.50){\v(3,-2){30.00}}
\end{picture}
\begin{picture}(114.00,90.00)
\put(0.00,24.00){$W_4:$}
\put(85.00,25.00){$\bu$}
\put(85.00,35.00){$0$}
\put(130.00,25.00){$\bu$}
\put(130.00,35.00){$4$}
\put(130.00,00.00){$\bu$}
\put(130.00,10.00){$5$}
\put(130.00,50.00){$\bu$}
\put(130.00,60.00){$3$}
\put(95.00,27.50){\v(1,0){30.00}}
\put(95.00,30.00){\v(3,2){30.00}}
\put(95.00,25.00){\v(3,-2){30.00}}
\put(40.00,0.00){$\bu$}
\put(40.00,10.00){$2$}
\put(40.00,50.00){$\bu$}
\put(40.00,60.00){$1$}
\put(50.00,5.50){\v(3,2){30.00}}
\put(50.00,49.50){\v(3,-2){30.00}}
\put(65.00,45.00){$a$}
\put(110.00,30.00){$d$}
\put(65.00,5.00){$b$}
\put(105.00,45.00){$c$}
\put(105.00,5.00){$e$}
\end{picture}

\begin{picture}(114.00,90.00)
\put(0.00,24.00){$W_5:$}
\put(85.00,25.00){$\bu$}
\put(85.00,35.00){$0$}
\put(65.00,45.00){$a$}
\put(65.00,5.00){$c$}
\put(55.00,30.00){$b$}
\put(105.00,45.00){$d$}
\put(105.00,5.00){$e$}
\put(130.00,00.00){$\bu$}
\put(130.00,10.00){$5$}
\put(130.00,50.00){$\bu$}
\put(130.00,60.00){$4$}
\put(95.00,30.00){\v(3,2){30.00}}
\put(95.00,25.00){\v(3,-2){30.00}}
\put(40.00,25.00){$\bu$}
\put(40.00,35.00){$2$}
\put(40.00,0.00){$\bu$}
\put(40.00,10.00){$3$}
\put(40.00,50.00){$\bu$}
\put(40.00,60.00){$1$}
\put(160.00,25.00){$cd=0$}
\put(50.00,27.50){\v(1,0){30.00}}
\put(50.00,5.50){\v(3,2){30.00}}
\put(50.00,49.50){\v(3,-2){30.00}}
\end{picture}

\begin{picture}(114.00,90.00)
\put(0.00,24.00){$W_6:$}
\put(40.00,0.00){$\bu$}
\put(40.00,10.00){$1$}
\put(90.00,0.00){$\bu$}
\put(90.00,10.00){$2$}
\put(140.00,0.00){$\bu$}
\put(140.00,10.00){$3$}
\put(190.00,0.00){$\bu$}
\put(190.00,-10.00){$4$}
\put(390.00,0.00){$\bu$}
\put(390.00,10.00){$8$}
\put(190.00,50.00){$\bu$}
\put(190.00,60.00){$9$}
\put(240.00,0.00){$\bu$}
\put(240.00,10.00){$5$}
\put(290.00,0.00){$\bu$}
\put(290.00,10.00){$6$}
\put(340.00,0.00){$\bu$}
\put(340.00,10.00){$7$}
\put(50.00,2.00){\v(1,0){35.00}}
\put(100.00,2.00){\v(1,0){35.00}}
\put(150.00,2.00){\v(1,0){35.00}}
\put(200.00,2.00){\v(1,0){35.00}}
\put(250.00,2.00){\v(1,0){35.00}}
\put(300.00,2.00){\v(1,0){35.00}}
\put(350.00,2.00){\v(1,0){35.00}}
\put(192.00,45.00){\v(0,-1){35.00}}
\put(62.00,-8.00){$a$}
\put(112.00,-8.00){$b$}
\put(162.00,-8.00){$c$}
\put(212.00,-8.00){$d$}
\put(262.00,-8.00){$e$}
\put(312.00,-8.00){$f$}
\put(362.00,-8.00){$g$}
\put(198.00,24.00){$h$}
\end{picture}

\bigskip
\bigskip

The wildness of the algebras $W_1-W_6$ follows from \cite{U}.

Let us consider the following box (see \cite{D} or \cite{Ro} for
definition) whose wildness  will be used in the proof of the
Theorem 1.

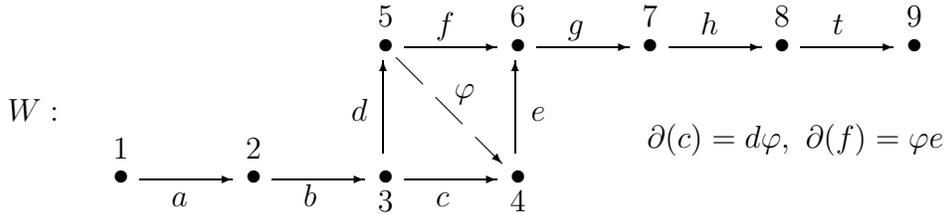
\begin{picture}(114.00,90.00)
\put(0.00,24.00){$W:$}
\put(40.00,0.00){$\bu$}
\put(40.00,10.00){$1$}
\put(90.00,0.00){$\bu$}
\put(90.00,10.00){$2$}
\put(140.00,0.00){$\bu$}
\put(140.00,-10.00){$3$}
\put(190.00,0.00){$\bu$}
\put(190.00,-10.00){$4$}
\put(140.00,50.00){$\bu$}
\put(140.00,60.00){$5$}
\put(190.00,50.00){$\bu$}
\put(190.00,60.00){$6$}
\put(240.00,50.00){$\bu$}
\put(240.00,60.00){$7$}
\put(290.00,50.00){$\bu$}
\put(290.00,60.00){$8$}
\put(340.00,50.00){$\bu$}
\put(340.00,60.00){$9$}
\put(50.00,2.00){\v(1,0){35.00}}
\put(150.00,2.00){\v(1,0){35.00}}
\put(100.00,2.00){\v(1,0){35.00}}
\put(150.00,52.00){\v(1,0){35.00}}
\put(200.00,52.00){\v(1,0){35.00}}
\put(250.00,52.00){\v(1,0){35.00}}
\put(300.00,52.00){\v(1,0){35.00}}
\put(142.00,12.00){\v(0,1){35.00}}
\put(192.00,12.00){\v(0,1){35.00}}
\put(147.00,48.00){\line(1,-1){10.00}}
\put(162.00,33.00){\line(1,-1){10.00}}
\put(177.00,18.00){\v(1,-1){10.00}}
\put(62.00,-8.00){$a$}
\put(112.00,-8.00){$b$}
\put(162.00,-8.00){$c$}
\put(162.00,57.00){$f$}
\put(212.00,57.00){$g$}
\put(262.00,57.00){$h$}
\put(312.00,57.00){$t$}
\put(130.00,24.00){$d$}
\put(198.00,24.00){$e$}
\put(169.00,33.00){$\f$}
\put(242.00,14.00){$\partial(c)=d\f ,\,\, \partial (f)=\f e$}
\end{picture}

\medskip
\bigskip
\bigskip

Consider the following
dimension vector $\vec d$: 

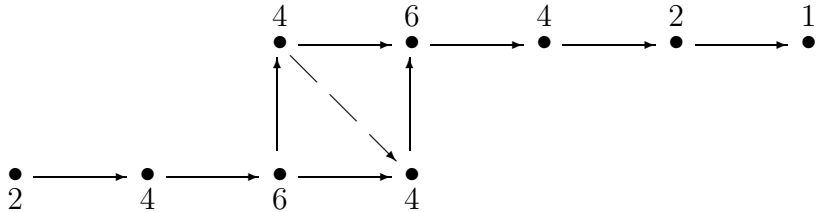
\begin{picture}(114.00,90.00)
\put(0.00,0.00){$\bu$}
\put(50.00,0.00){$\bu$}
\put(100.00,0.00){$\bu$}
\put(150.00,0.00){$\bu$}
\put(100.00,50.00){$\bu$}
\put(150.00,50.00){$\bu$}
\put(200.00,50.00){$\bu$}
\put(250.00,50.00){$\bu$}
\put(300.00,50.00){$\bu$}
\put(10.00,2.00){\v(1,0){35.00}}
\put(110.00,2.00){\v(1,0){35.00}}
\put(60.00,2.00){\v(1,0){35.00}}
\put(110.00,52.00){\v(1,0){35.00}}
\put(160.00,52.00){\v(1,0){35.00}}
\put(210.00,52.00){\v(1,0){35.00}}
\put(260.00,52.00){\v(1,0){35.00}}
\put(102.00,12.00){\v(0,1){35.00}}
\put(152.00,12.00){\v(0,1){35.00}}
\put(107.00,48.00){\line(1,-1){10.00}}
\put(122.00,33.00){\line(1,-1){10.00}}
\put(137.00,18.00){\v(1,-1){10.00}}
\put(0.00,-10.00){$2$}
\put(50.00,-10.00){$4$}
\put(100.00,-10.00){$6$}
\put(150.00,-10.00){$4$}
\put(100.00,59.00){$4$}
\put(150.00,59.00){$6$}
\put(200.00,59.00){$4$}
\put(250.00,59.00){$2$}
\put(300.00,59.00){$1$}
\end{picture}

\bigskip
\bigskip

Since $f_{W}(\vec d)=-1$, the wildness of the box $W$ follows from \cite{D}.

\bigskip
\bigskip

\subsection{}

\begin{theorem}
\begin{itemize}
\item[(i)]
 Let $S=S(n,d)$ be a Schur algebra. Then $S$ is derived tame if and only if 
every block of $S$ is Morita equivalent to $A_1$ or $A_2$;

\item[(ii)] Let $S=S(n,d)_r$ be an infinitesimal Schur algebra.
  Then $S$ is derived tame if and only if every block of $S$ is Morita equivalent to $A_1$ or $F_3$.
\end{itemize}
\end{theorem}

\begin{proof}
 It follows from \cite{DEMN} that if $S(n,d)$ or $S(n,d)_r$ is not wild then any of 
 its non-semisimple blocks 
is Morita equivalent 
to one of the algebras in Section 2.
We show that all algebras in this list with the exception of $A_1$, $A_2$
and $F_3$ are derived wild. Let $A$ be one of the algebras from Section 2 and let 
  $B$ be one of the algebras $W_1-W_6$. 
Since $B$ is wild, there exists $B-k\left< x,y\right>$-bimodule $M=M(B)$
such that the functor $M\otimes_{k\left< x,y\right> }-$ preserves
indecomposability and isomorphism classes. 
We set $d_i=\dim M(i)$ and  denote by
$[M(x)]$  the matrix corresponding to the map
$M(x): M(s(x))\rightarrow M(e(x))$
 with respect to some fixed basis.

\medskip
\noindent (1) Let $A=G$,  $M=M(W_1)$
and let  $N^{\bullet}$ be the
following complex of $A-k\left< x,y\right>$-bimodules:

$$
\cdots\rightarrow 0\rightarrow N^{1}=\oplus_{i=1}^{3}{\bf P}_i^{d_i}
\stackrel{\partial^{1}}{\rightarrow}{\bf P}_0^{d_0}\stackrel{\partial^{2}}{\rightarrow}
{\bf P}_{1}^{d_4}\oplus {\bf P}_{2}^{d_5}\rightarrow 0\rightarrow \cdots,
$$

\noindent where 
$$\partial^{1}=
\left( \begin{array}{lll}
p(\a_1)[M(a)]&p(\a_2)[M(b)]&p(\a_3)[M(c)]
\end{array}
\right)^{T}, 
\partial^{2}=
\left( \begin{array}{ll}
p(\b_1)[M(d)]&p(\b_2)[M(e)]
\end{array}
\right).$$

\medskip
\noindent (2) Let $A=B$, 
 $M=M(W_2)$ and
 let $N^{\bullet}$ be the
following complex of $A-k\left<x,y\right>$-bimodules:

$$
\cdots\rightarrow 0\rightarrow N^{1}=\oplus_{i=1}^{4}{\bf P}_i^{d_i}
\stackrel{\partial^{1}}{\rightarrow} {\bf P}_0^{d_0}\stackrel{\partial^{2}}{\rightarrow}
{\bf P}_1^{d_5}\rightarrow 0\rightarrow \cdots,
$$

\noindent where 
$$\partial^{1}=
\left( \begin{array}{llll}
p(\a_1)[M(a)]&p(\a_2)[M(b)]&p(\a_3)[M(c)]&p(\a_4)[M(d)]
\end{array}
\right)^{T}, \partial^{2}=p(\b_1)[M(e)].$$

\medskip
\noindent (3) Let $A=B_1$,  
 $M=M(W_3)$ and let
  $N^{\bullet}$ be the
following complex of $A-k\left< x,y\right>$-bimodules:

$$
\cdots\rightarrow 0\rightarrow N^{1}=\oplus_{i=1}^{4}{\bf P}_i^{d_i}
\stackrel{\partial^{1}}{\rightarrow} {\bf P}_0^{d_0}\stackrel{\partial^{2}}{\rightarrow}
{\bf P}_1^{d_5}\rightarrow 0\rightarrow \cdots,
$$

\noindent where 
$$\partial^{1}=
\left( \begin{array}{llll}
p(\a_1)[M(a)]&p(\a_2)[M(b)]&p(\a_3)[M(c)]&p(\a_4)[M(d)]
\end{array}
\right)^{T}, \partial^{2}=p(\b_4)[M(e)].$$

\medskip
\noindent (4) Let $A=D$, 
 $M=M(W_4)$ and let
  $N^{\bullet}$ be the
following complex of $A-k\left< x,y\right>$-bimodules:

$$
\cdots\rightarrow 0\rightarrow N^{1}={\bf P}_1^{d_1}\oplus {\bf P}_2^{d_2}
\stackrel{\partial^{1}}{\rightarrow} {\bf P}_0^{d_0}\stackrel{\partial^{2}}{\rightarrow}
{\bf P}_1^{d_3}\oplus {\bf P}_2^{d_4}\oplus {\bf P}_0^{d_5}\rightarrow 0\rightarrow \cdots,
$$

\noindent where 
$$\partial^{1}=
\left( \begin{array}{ll}
p(\a_1)[M(a)]&p(\a_2)[M(b)]
\end{array}
\right)^{T},
\partial^{2}=
\left( \begin{array}{lll}
p(\b_1)[M(c)]&p(\b_2)[M(d)]&p(\b_3\a_3)[M(e)]
\end{array}
\right).$$

\medskip
\noindent (5) Let $A=D_3$, 
 $M=M(W_6)$ and let
  $N^{\bullet}$ be the
following complex of $A-k\left< x,y\right>$-bimodules:

Set $N^{i}={\bf P}_2^{d_i}$ for $i\in \{1,2,4,5,6,7,8\}$,
$N^{3}={\bf P}_2^{d_3}\oplus {\bf P}_1^{d_9}$,
$\partial^{1}=p(\b_2\a_2)[M(a)]$,
$\partial^{2}=
\left( \begin{array}{ll}
p(\b_2\a_2)[M(b)]&0
\end{array}
\right)$,
$\partial^{3}=
\left( \begin{array}{ll}
p(\b_2\a_2)[M(c)]&p(\a_1)[M(h)]
\end{array}
\right)^{T}$,
$\partial^{4}=p(\b_2\a_2)[M(d)]$,
$\partial^{5}=p(\b_2\a_2)[M(e)]$,
$\partial^{6}=p(\b_2\a_2)[M(f)]$ and
$\partial^{7}=p(\b_2\a_2)[M(g)]$.

\medskip
\noindent (6) Let $A=D_4$,  
 $M=M(W_5)$ and let
  $N^{\bullet}$ be the
following complex of $A-k\left< x,y\right>$-bimodules:

$$
\cdots\rightarrow 0\rightarrow N^{1}=\oplus_{i=1}^{3}{\bf P}_i^{d_i}
\stackrel{\partial^{1}}{\rightarrow} {\bf P}_0^{d_0}\stackrel{\partial^{2}}{\rightarrow}
{\bf P}_1^{d_4}\oplus {\bf P}_2^{d_5}\rightarrow 0\rightarrow \cdots,
$$

\noindent where 
$$\partial^{1}=
\left( \begin{array}{lll}
p(\a_1)[M(a)]&p(\a_3)[M(b)]&p(\a_2)[M(c)]
\end{array}
\right)^{T},
\partial^{2}=
\left( \begin{array}{lll}
p(\b_1)[M(d)]&p(\b_2\a_2)[M(e)]
\end{array}
\right).$$

\medskip
\noindent (7) Let $A=H_4$, 
 $M=M(W_5)$ and let 
  $N^{\bullet}$ be the
following complex of $A-k\left< x,y\right>$-bimodules:

$$
\cdots\rightarrow 0\rightarrow N^{1}={\bf P}_1^{d_1}\oplus {\bf P}_3^{d_2}
\oplus {\bf P}_2^{d_3}
\stackrel{\partial^{1}}{\rightarrow} {\bf P}_0^{d_0}\stackrel{\partial^{2}}{\rightarrow}
{\bf P}_1^{d_4}\oplus {\bf P}_2^{d_5}\rightarrow 0\rightarrow \cdots,
$$

\noindent where 
$$\partial^{1}=
\left( \begin{array}{lll}
p(\a_1)[M(a)]&p(\a_3)[M(b)]&p(\a_2)[M(c)]
\end{array}
\right)^{T},
\partial^{2}=
\left( \begin{array}{ll}
p(\b_3)[M(d)]&p(\b_1)[M(e)]
\end{array}
\right).$$

\medskip
\noindent (8) Let $A\in \{A_m, m>2, F_r, r>3, R_4\}$. 

Since box $W$ is wild, there exists $W-k\left< x,y\right>$-bimodule $M$
such that the functor $M\otimes_{k\left< x,y\right> }-$ preserves
indecomposability and isomorphism classes.
 Denote by $N^{\bullet}$ the
following complex of $A-k\left< x,y\right>$-bimodules.

Set $N^{i}={\bf P}_3^{d_i}$ for $1\leq i\leq 3$,
$N^{j}={\bf P}_2^{d_{j+1}}$ for $5\leq j\leq 8$,
$N^{4}={\bf P}_2^{d_5}\oplus {\bf P}_3^{d_4}$,
$$\partial^{1}=p(\b_2\a_2)[M(a)],
\partial^{2}=p(\b_2\a_2)[M(b)],
\partial^{3}=
\left( \begin{array}{ll}
p(\b_2)[M(d)]&p(\b_2\a_2)[M(c)]
\end{array}
\right),$$
$$\partial^{4}=
\left( \begin{array}{ll}
p(\a_2\b_2)[M(f)]&p(\b_2)[M(e)]
\end{array}
\right)^{T},
\partial^{5}=p(\a_2\b_2)[M(g)],$$
$$\partial^{6}=p(\a_2\b_2)[M(h)],
\partial^{7}=p(\a_2\b_2)[M(t)].$$

\bigskip
\medskip
It is not difficult to verify that the functor 
$N^{\bullet}\otimes_{k\left< x,y\right> }-$, which acts from the category of
finite-dimensional $k\left< x,y\right>$-modules to the category ${\frak p}(A)$,
where $A$ is one of the algebras from (1)-(8),
preserves indecomposability and isomorphism classes.
So, $A$ is derived wild.

It follows from \cite{H} that algebra $A_1$ is derived tame.
It follows from \cite{BM} that the algebra $A_2$ is derived tame.
The derived tameness of the algebra $F_3$ follows from Theorem 4
(see Section 4).

\end{proof}

\subsection{\bf Proof of Theorem 1.}

It follows from Theorem 2 that a Schur (resp., infinitesimal Schur) 
algebra $S(n,d)$ (resp., $S(n,d)_r$) is derived tame if all its blocks are of type $A_1$ 
or $A_2$ (resp., $A_1$ or $F_3$). Hence, in order to classify derived tame Schur 
(resp., infinitesimal Schur)  algebras 
it is enough to choose among representation tame and representation finite algebras 
those that have all blocks of type $A_1$, $A_2$ or $F_3$.   

Let $S(n,d)$ (resp., $S(n,d)_r$) be a Schur (resp., infinitesimal Schur) algebra. 
It follows from \cite{DN} that it is semisimple if and only if it satisfies conditions  
 a), b) in (i) (resp., a), b), c) in (ii)). If $S(n,d)$ is of type c) in (i) then any of 
its non-semisimple blocks is Morita equivalent to $A_2$ by 5.4, 5.5 and 1.3 in \cite{E}. Suppose now that $S(n,d)$ is a Schur algebra 
of finite representation type with $n=2$ and $d<p^2$ \cite{E}. Such algebra has all non-semisimple blocks of type $A_2$ if and only if it satisfies d) in (i) by Proposition 5.1 in \cite{E}. Let $S(n,d)$ be a Schur 
algebra of finite representation type with $n\geq 3$ and $d<2p$. The basic algebra of $S(n,d)$ is a direct sum 
of all blocks of the group algebra $kS_d$ of the symmetric group $S_d$ (see 1.4 in \cite{E}). The irreducible representations of $S_d$ are parametrized by the partitions of $d$. Moreover, two such representations belong to the same block if the corresponding partitions have the same $p$-core \cite{JK}. By  4.1 in \cite{E}, if a block 
of $kS_d$ has $s$ partitions with $\leq n$ parts then it is equivalent to the algebra $A_s$. Applying this 
 we conclude that $S(n,d)$ with $n\geq 3$ and $d<2p$ has all non-semisimple blocks Morita equivalent to $A_2$ 
if and only if it satisfies the conditions e) and f) in (i). 

Now let $S(n,d)_r$ be an infinitesimal Schur algebra satisfying 
d) in (ii). Then its non-semisimple blocks are Morita equivalent 
to $F_3$ by \cite{DNP2}, Section 6.3. Finally, the remaining infinitesimal Schur algebras 
of finite representation type have non-semisimple 
blocks Morita equivalent to the algebra $F_3$ 
if and only if they satisfy e) and f) in (ii) by Propositions 2.4 and 3.2 in 
\cite{DNP3}. This completes the proof of Theorem 1.

\section{Indecomposables in derived categories of algebras $A_1$, $A_2$ and 
$F_3$}

\subsection{$A=A_1$.}

It follows from \cite{H} that
the projective complexes 
$$P_i^{\bullet}:\,\,\,\,\,\, \cdots\rightarrow 0\rightarrow P^{i}={\bf P}_1
\rightarrow 0\rightarrow \cdots ,$$

\noindent where $i\in \Z$, 
constitute an exhaustive list of pairwise non-isomorphic
indecomposable objects of ${\bf D}^{b}(A_1)$.

\begin{corollary}
An algebra of type $A_1$ is derived finite.
\end{corollary}

\subsection{$A=A_2$.} 

Let $\a=\a_1, \b=\b_1$ and let $e_1, e_2$ be the idempotents 
corresponding to the vertices.
Consider the following projective complexes of ${\bf D}^{b}(A_2)$.

\begin{itemize}
\item 
$P^{\bullet}_{e_1}:\,\,\,\,\,\, \cdots\rightarrow 0\rightarrow P^{1}={\bf P}_1
\rightarrow 0\rightarrow \cdots$
\item 
$P^{\bullet}_{e_2}:\,\,\,\,\,\, \cdots\rightarrow 0\rightarrow P^{1}={\bf P}_2
\rightarrow 0\rightarrow \cdots$
\item 
$P^{\bullet}_{\a}:\,\,\,\,\,\, \cdots\rightarrow 0\rightarrow P^{1}={\bf P}_1
\stackrel{p(\a)}{\rightarrow} {\bf P}_2\rightarrow 0\rightarrow \cdots$
\item 
$P^{\bullet}_{\b}:\,\,\,\,\,\, \cdots\rightarrow 0\rightarrow P^{1}={\bf P}_2
\stackrel{p(\b)}{\rightarrow} {\bf P}_1\rightarrow 0\rightarrow \cdots$
\item 
$P^{\bullet}_{(\b\a)^{s}}:\,\,\,\,\,\, \cdots\rightarrow 0\rightarrow P^{1}={\bf P}_2
\stackrel{p(\b\a)}{\rightarrow} {\bf P}_2\stackrel{p(\b\a)}{\rightarrow}  \cdots
\stackrel{p(\b\a)}{\rightarrow} P^{s+1}={\bf P}_2\rightarrow 0\rightarrow \cdots$
\item 
$P^{\bullet}_{\a(\b\a)^{s}}:\,\,\,\,\,\, \cdots\rightarrow 0\rightarrow
 P^{1}={\bf P}_1
\stackrel{p(\a)}{\rightarrow}
{\bf P}_2
\stackrel{p(\b\a)}{\rightarrow} {\bf P}_2\stackrel{p(\b\a)}{\rightarrow}  \cdots
\stackrel{p(\b\a)}{\rightarrow} P^{s+2}={\bf P}_2\rightarrow 0\rightarrow \cdots$
\item 
$P^{\bullet}_{(\b\a)^{s}\b}:\,\,\,\,\,\, \cdots\rightarrow 0\rightarrow
 P^{1}={\bf P}_2
\stackrel{p(\b\a)}{\rightarrow} {\bf P}_2\stackrel{p(\b\a)}{\rightarrow} \cdots
\stackrel{p(\b\a)}{\rightarrow} {\bf P}_2
\stackrel{p(\b)}{\rightarrow} P^{s+2}={\bf P}_1
\rightarrow 0\rightarrow \cdots$
\item 
$P^{\bullet}_{\a(\b\a)^{s}\b}:\,\,\,\,\,\, \cdots\rightarrow 0\rightarrow
 P^{1}={\bf P}_1
\stackrel{p(\a)}{\rightarrow}
{\bf P}_2
\stackrel{p(\b\a)}{\rightarrow} {\bf P}_2\stackrel{p(\b\a)}{\rightarrow}  \cdots \stackrel{p(\b\a)}{\rightarrow} 
{\bf P}_2\stackrel{p(\b\a)}{\rightarrow} {\bf P}_2
\stackrel{p(\b)}{\rightarrow} P^{s+3}={\bf P}_1
\rightarrow 0\rightarrow \cdots$
\end{itemize}

\begin{theorem}
The projective complexes $T^{i}(P^{\bullet}_{e_1})$,
$T^{i}(P^{\bullet}_{e_2})$,
$T^{i}(P^{\bullet}_{\a})$,
$T^{i}(P^{\bullet}_{\b})$,
$T^{i}(P^{\bullet}_{(\b\a)^{s}})$,\\
$T^{i}(P^{\bullet}_{\a(\b\a)^{s}})$,
$T^{i}(P^{\bullet}_{(\b\a)^{s}\b})$ and
$T^{i}(P^{\bullet}_{\a(\b\a)^{s}\b})$,
where $i\in \Z$, $s\in \N$,
 constitute an exhaustive list of   pairwise non-isomorphic
indecomposable objects of ${\bf D}^{b}(A_2)$.
\end{theorem}
\begin{proof}
The proof follows from \cite{BM}.
\end{proof}

\begin{corollary}
An algebra of type $A_2$ is derived discrete.
\end{corollary}

\subsection{$A=F_3$.}

We will construct a certain box (see \cite{D} or \cite{Ro} for
definition) corresponding to the algebra $F_3$.
Consider the path algebra $B=kQ/I$, where
$Q_0=\{1[i], 2[i], 3[i]\, |\,$ 
$i\in \Z\}$,
$Q_1=\{a[i], b[i], c[i], d[i], f[i]\, |\,$
$i\in \Z\}$,
$s(a[i])=s(c[i])=s(f[i])=e(b[i-1])=e(d[i-1])=e(f[i-1])=2[i]$,
$s(b[i])=e(a[i-1])=1[i]$,
$s(d[i])=e(c[i-1])=3[i]$ and \\
$I=<a[i]b[i+1]+c[i]d[i+1]\, |\,$
$i\in \Z>$. 

Consider a normal box ${\cal B}=(B,V)$
with a kernel $\overline{V}$ freely generated by the set
$\{ \a[i], \b[i], \g[i], \d[i], \f[i]\, |\,$
$i\in\Z\}$, where
$s(\a[i])=s(\g[i])=s(\f[i])=e(\b[i])=e(\d[i])=e(\f[i])=2[i]$,
$s(\b[i])=e(\a[i])=1[i]$,
$s(\d[i])=e(\g[i])=3[i]$ and
with the differential $\partial$ given by the formulas:
$$\partial (f[i])=a[i]\b[i+1]+c[i]\d[i+1]+\a[i]b[i]+\g[i]d[i],$$
$$\partial (\f[i])=\a[i]\b[i]+\g[i]\d[i],
\partial (x)=0 \,\,{\rm for}\,\, x\not\in \{f[i],\f[i]\,|\,i\in\Z\}.$$

The category of the representations of the box ${\cal B}$ will be
denoted by ${\rm rep}\,({\cal B})$.

\begin{proposition}
 The category $\frak p (F_3)$ is equivalent to ${\rm rep}\,({\cal B})$.
\end{proposition}
\begin{proof}
 Consider the functor 
$G: {\rm rep}\, ({\cal B})\rightarrow \frak p(F_3)$
defined as follows.
The modules are
$$G(M)^{j}=\oplus_{i\in \{1,2,3\}}{\bf P}_i^{\dim_k M(i[j])}$$ 
and the differential maps are
$$\partial^{j}_M=
\left( \begin{array}{lll}
0&p(\a_1)[M(b[j])]&0\\
p(\b_1)[M(a[j])]&p(\b_1\a_1)[M(f[j])]&p(\a_2)[M(c[j])]\\
0&p(\b_2)[M(d[j])]&0
\end{array}
\right),$$
\noindent where $[M(x[j])]$ is the matrix corresponding to the map
$M(x[j])$ in some fixed basis.
It is easy to see that $G$ is a representation equivalence.
\end{proof}

We recall some definitions and results related to the bunches
of semi-chains considered by Bondarenko in \cite{B} and Deng
in \cite{De}
in a form convenient for our purposes
(see also \cite{CB} for an alternative approach).
We will use the classification of indecomposables  representations of a
bunch of semi-chains given in \cite{B}.
We will use some notation from \cite{DG}.

\begin{definition}
A {\it bunch of semi-chains} ${\bf C}=\{{\bf I}, E_i,F_i,\sim\}$ is defined
by following data:
\begin{enumerate}
\item A set ${\bf I}$ of indices;

\item Two semi-chains (i.e., partially ordered sets with the condition 
that  each element is incomparable with at most one other) 
$E_i$ and $F_i$ given for each $i\in {\bf I}$;

Put ${\bf E}:=\cup_{i\in {\bf I}}E_i$,
${\bf F}:=\cup_{i\in {\bf I}}F_i$ and $|{\bf C}|:={\bf E}\cup {\bf F}$.

\item An equivalence relation $\sim$ on 
$|{\bf C}|$ such that each equivalence class consists of at most
2 elements and if $a\sim b\not=a$ for some $a\in E_i$ 
(resp., $a\in F_i$),  then $a$  is
comparable with all elements of $E_i$ (resp., $F_i$).
\end{enumerate}
\end{definition}

We consider the ordering on $|{\bf C}|$, which is just
the union of all orderings on $E_i$ and $F_i$
(i.e., $a<b$ means that $a$ and $b$ belong to the same
semi-chain $E_i$ or $F_i$ and $a<b$ in this semi-chain).

We  construct a certain box associated with a  
bunch of semi-chains ${\bf C}$.
Given $u\in |{\bf C}|$ we denote by $\bar u$ the corresponding element 
of $|{\bf C}|/\sim$.
Consider the path algebra $A({\bf C})=kQ=kQ({\bf C})$, where
$Q_0=|{\bf C}|/\sim,$ 
$Q_1=\{a_u^{v}\, |\,$
$u\in E_i, v\in F_i, i\in \Z\}$,
(we also will assume that $a_u^{v}=0$ for
all other cases),
$s(a_u^{v})=\bar u$ and
$e(a_u^{v})=\bar v$. 

Consider a normal box ${\cal C}=(A({\bf C}),V=V({\bf C}))$
with a kernel $\overline{V}$ freely generated by the set
$\{ \f_u^{v}, \psi_m^{n}\,|\,$
$u, v\in E_i, m, n\in F_i, u<v, n<m, i\in\Z\}$
(we also will suppose that $\f_u^{v}=\psi_m^{n}=0$ for
all other cases), where
$s(\f_u^{v})=\bar u$,
$s(\psi_m^{n})=\bar m$,
$e(\f_u^{v})=\bar v$,
$e(\psi_m^{n})=\bar n$  and
with the differential $\partial$ given by the formulas:
$$\partial (a_u^{v})=
\sum_{n>u}\f_u^{n}a_n^{v}
+\sum_{m>v}a_u^{m}\psi_m^{v},$$
$$\partial (\f_u^{v})=\pm\sum_{u<n<v}\f_u^{n}\f_n^{v},$$
$$\partial (\psi_u^{v})=\pm\sum_{v<n<u}\psi_u^{n}\psi_n^{v}.$$

The choice of signs in the last formulas guarantees the condition
$\partial^{2}=0$.

The category of representations of the bunch of semi-chains
${\bf C}$ is then 
defined as the category ${\rm rep}\,{\cal C}$.
One can easily verify that this definition gives just the same
representations as the definition in  \cite{B}.

Consider an equivalence relation $\sim_c$ on $|{\bf C}|$
given by the following rule:
$a\sim_c b$ if and only if either $a=b$ or 
$a$ and $b$ belong to the same semi-chain $E_i$ or $F_i$ and 
$a$ is incomparable with $b$. 
In case the equivalence class $x\in |{\bf C}|/\sim_c$ consists
of a unique element $a\in |{\bf C}|$ we will identify $x$  with $a$
and write $x\in |{\bf C}|$. Consider a relation  
$\sim$ on $|{\bf C}|/\sim_c$ given by the following rule: 
$a\sim b$ if and only if either $a=b$ and $a$ consists of two elements or 
$a,b\in |{\bf C}|$ and $a\sim b\ne a$.    

\begin{definition} Let ${\bf C}=\{{\bf I}, E_i,F_i,\sim\}$ be
a bunch of semi-chains.
\begin{itemize}
\item  A {\bf C}-{\rm word} is a sequence $w=w_0r_1w_1r_2w_2\cdots r_mw_m$,
where $w_k\in |{\bf C}|/\sim_c$ and each $r_k$ is either $\sim$ or $-$,
such that  for all possible values of $k$:

(a) $r_k= \sim$ if and only if  $w_{k-1}\sim w_k$;

(b) $r_k=-$ if and only if $w_{k-1}\subset E_i$, $w_{k}\subset F_i$ 
for some $i\in \Z$ or vice versa;

(c) $r_{k+1}\ne r_k$.

Possibly $m=0$, i.e., $w\in |{\bf C}|/\sim_c$. 

\item Call a {\bf C}-word  $w=w_0r_1w_1r_2w_2\cdots r_mw_m$ 
a {\bf C}-{\rm cycle}  if $w_m=w_0$, $r_1= \sim$ and $r_m=-$.

\item Call a {\bf C}-word  $w=w_0r_1w_1r_2w_2\cdots r_mw_m$
{\rm full} if, whenever $w_0\in |{\bf C}|$ and $w_0$ 
is not a unique element in its
equivalence class $\overline{w_0}$, 
then $r_1= \sim$, and whenever $w_m\in |{\bf C}|$ and $w_m$ is not
a  unique element in its equivalence class $\overline{w_m}$, 
then $r_m= \sim$.

\item Call a {\bf C}-cycle  $w=w_0r_1w_1r_2w_2\cdots r_mw_m$
{\rm aperiodic} if the sequence  $w_0r_1w_1\cdots r_m$
can not be written as  a multiple self-concatenation $v\cdots v$
of a shorter sequence $v$.

Given a {\bf C}-word  $w=w_0r_1w_1r_2w_2\cdots r_mw_m$ 
set $w^{\ast}:=w_mr_m\cdots r_2w_1r_1w_0$.
 
\item Call a {\bf C}-word  $w=w_0r_1w_1r_2w_2\cdots r_mw_m$ 
{\rm simple} if, from  $w=u\sim u^{\ast}\sim u\cdots$ for some
{\bf C}-word $u$, it follows that $w=u$.

\end{itemize}
\end{definition}

Denote  by
$Ind \,k[x]$ the set of indecomposable polynomials with highest coefficient
$1$ except  $\{x^d |d\geq 1\}$.

\begin{definition}  Let ${\bf C}=\{{\bf I}, E_i,F_i,\sim\}$ be
a bunch of semi-chains.

\begin{itemize}
\item Given a full {\bf C}-word $w=w_0r_1w_1r_2w_2\cdots r_mw_m$,
we set $d_l(w)=1$ if  $w_0\sim w_0$ and either
$r_1=-$ or $m=0$, and we set $d_l(w)=0$ otherwise. 

\item Given a full {\bf C}-word $w=w_0r_1w_1r_2w_2\cdots r_mw_m$,
we set $d_r(w)=1$ if  $w_m\sim w_m$, $m>0$ and 
$r_m=-$, and we set $d_r(w)=0$ otherwise. 

\item By a {\rm usual string} we mean a simple full {\bf C}-word 
$w=w_0r_1w_1r_2w_2\cdots r_mw_m$
such that $d_l(w)+d_r(w)=0$.

\item By a {\rm special string} we mean a pair $(w,k)$, where
 $w=w_0r_1w_1r_2w_2\cdots r_mw_m$ 
is a simple full {\bf C}-word 
such that $0<d_l(w)+d_r(w)<2$ and $k\in \{0,1\}$.

\item By a {\rm bispecial string} we mean a quadruple $(w,k,l,n)$, where
$w=w_0r_1w_1r_2w_2\cdots r_mw_m$ 
is a simple full {\bf C}-word 
such that $d_l(w)+d_r(w)=2$ and $k,l\in \{0,1\}$, $n\in \N$.

\item By a {\rm string} we mean usual or special or bispecial string. 

\item By a {\rm band} we mean a pair $(w,f)$, where
$f\in Ind\, k[x]$ and $w$ is an aperiodic {\bf C}-cycle.
\end{itemize}
\end{definition}

For each string and each band,  
Bondarenko constructed in \cite{B} some
indecomposable representation of  {\bf C}
and proved that up to some isomorphisms between of them
(see \cite{B} for details) they constitute an exhaustive list of pairwise non-isomorphic
indecomposable representations.

>From now on we will consider the bunch of semi-chains
${\bf C}={\bf C}(F_3):=\{\Z, E_i,F_i, \sim\}$ , where
$E_i=\{y[i]^{-}<x[i]<z[i], y[i]^{+}<x[i]\}$,
$F_i=\{r[i]<p[i]^{-}<q[i]>p[i]^{+}>r[i]\}$
and the equivalence relation $\sim$  is given by 
$r[i]\sim x[i+2]$, $q[i]\sim z[i+1]$.
We denote by $p[i]$ (resp., $y[i]$) the class 
$\{p[i]^{-}, p[i]^{+}\}\in |{\bf C}(F_3)|/\sim_c$
(resp., $\{y[i]^{-}, y[i]^{+}\}\in |{\bf C}(F_3)|/\sim_c$).
 We denote by ${\cal C}(F_3)$ the corresponding box. 
 
We associate to indecomposable representations of 
${\bf C}(F_3)$ certain finite projective complexes which, as we shall see,
give all indecomposables in the category $\frak p (F_3)$.

\begin{definition}
Let $M$ be an indecomposable representation of the bunch of semi-chains
${\bf C}(F_3)$. Then $P(M)^{\bullet}$ is the projective complex
$\cdots\rightarrow P(M)^{i}\stackrel{\partial_{P(M)^{\bullet}}^{i}}{\rightarrow}P(M)^{i+1}
\rightarrow \cdots$
defined as follows.
The modules are
$$P(M)^{i}={\bf P}_1^{d_{1,i}}\oplus {\bf P}_2^{d_{2,i}}
\oplus{\bf P}_3^{d_{3,i}}$$ and
the differential maps are
$$\partial_{P(M)^{\bullet}}^{i}=
\left( \begin{array}{lll}
0&p(\a_1)A_i&0\\
p(\b_1)B_i&p(\b_1\a_1)C_i&p(\a_2)D_i\\
0&p(\b_2)R_i&0
\end{array}
\right),$$
where
$$A_i=
\left( \begin{array}{llll}
M_{x[i]}^{r[i]}&M_{x[i]}^{p[i]^{-}}&M_{x[i]}^{p[i]^{+}}&M_{x[i]}^{q[i]}\\
0&0&0&0\\
M_{y[i]^{+}}^{r[i]}&M_{y[i]^{+}}^{p[i]^{-}}&M_{y[i]^{+}}^{p[i]^{+}}&M_{y[i]^{+}}^{q[i]}
\end{array}
\right),
R_i=
\left( \begin{array}{llll}
-M_{x[i]}^{r[i]}&-M_{x[i]}^{p[i]^{-}}&-M_{x[i]}^{p[i]^{+}}&-M_{x[i]}^{q[i]}\\
0&0&0&0\\
M_{y[i]^{-}}^{r[i]}&M_{y[i]^{-}}^{p[i]^{-}}&M_{y[i]^{-}}^{p[i]^{+}}&M_{y[i]^{-}}^{q[i]}
\end{array}
\right),$$
$$C_i=
\left( \begin{array}{llll}
0&0&0&0\\
0&0&0&0\\
0&0&0&0\\
M_{z[i]}^{r[i]}&M_{z[i]}^{p[i]^{-}}&M_{z[i]}^{p[i]^{+}}&M_{z[i]}^{q[i]}
\end{array}
\right),
B_i=
\left( \begin{array}{lll}
E&0&0\\
0&E&0\\
0&0&0\\
0&0&0
\end{array}
\right),
D_i=
\left( \begin{array}{lll}
E&0&0\\
0&0&0\\
0&E&0\\
0&0&0
\end{array}
\right),$$
\noindent $M_{u[i]}^{v[i]}$ is the matrix corresponding to the map
$M(a_{u[i]}^{v[i]})$ in some fixed basis
and where $d_{i,j}$ and dimensions of all blocks are uniquely
defined by dimensions of  the matrix $M_{u[i]}^{v[i]}$.
\end{definition}

\begin{proposition}
The projective complexes $P(M)^{\bullet}$,
where $M\in {\bf Ver}\, {\rm rep}\,({\cal C}(F_3))$, 
constitute an exhaustive list of pairwise non-isomorphic
indecomposable objects of $\frak p (F_3)$.
\end{proposition}

\begin{proof}
Let ${\cal A}=(A,V)$ be the box corresponding to $F_3$ (see above).
Let us consider the sub-box ${\cal A}^{\prime}=(A^{\prime},V^{\prime})$
of ${\cal A}$, where $A^{\prime}_0=A_0$, 
$A^{\prime}_1=\{a[i],c[i]\,|\,i\in\Z\}$ and $\overline{V^{\prime}}=0$.
It is easy to see that 
${\cal A}^{\prime}\cong\amalg_{\Z}{\cal A}^{\prime\prime}$,
where ${\cal A}^{\prime\prime}$ is the principal box corresponding to the hereditary
algebra $A^{\prime\prime}:\, \bullet\leftarrow\bullet\rightarrow\bullet$
(i.e., ${\cal A}^{\prime\prime}=(A^{\prime\prime},A^{\prime\prime})$).

Consider the trivial category $D^{\prime}$ 
with the set of  vertices  $|{\bf C}(F_3)|/_{\sim}$ and
the functor $F^{\prime}: A^{\prime}\rightarrow D^{\prime}$ which maps:
$$
1[i]\rightarrow \overline{x[i+1]}\oplus \overline{p[i-1]^{+}}
\oplus \overline{y[i+1]^{+}},
$$
$$
2[i]\rightarrow \overline{x[i+1]}\oplus \overline{p[i-1]^{+}}
\oplus \overline{p[i-1]^{-}}
\oplus \overline{z[i]},
$$
$$
3[i]\rightarrow \overline{x[i+1]}\oplus \overline{p[i-1]^{-}}
\oplus \overline{y[i+1]^{-}},
$$
$$
a[i]\rightarrow  \left( \begin{array}{lll}
1&0&0\\
0&1&0\\
0&0&0\\
0&0&0
\end{array}
\right),\,\,\,\,
c[i]\rightarrow  \left( \begin{array}{lll}
1&0&0\\
0&0&0\\
0&1&0\\
0&0&0
\end{array}
\right).
$$

Construct the
{amalgamation} $B=A\amalg^{A^{\prime}}B^{\prime}\,$, or, the same, the
couniversal square:
$$
\begin{array}{ccc}
     A^{\prime}     & \stackrel{F^{\prime}}{\longrightarrow} &       B^{\prime}       \\
        \downarrow   &               &      \downarrow   \\
     A      & \stackrel{F}{\longrightarrow}   &       B       
\end{array}
$$
Consider now the
box ${\cal A}^F=(B,V^{F})$, where
$V^{F}={\cal B}\otimes_{\cal A}V\otimes_{\cal A}{\cal B}$.

We say that the box ${\cal A}^{F}$ is obtained from ${\cal A}$
{\em by reducing the sub-box} ${\cal A}^{\prime}$. 
By \cite{D}, $F$ induces  the representation equivalence
$F^{\ast}: {\rm rep}\,({\cal A}^{F})\rightarrow {\rm rep}\,({\cal A})$.

Straightforward calculation, which we omit, shows that
$B=A({\bf C}(F_3))$ and
there exists the sub-bimodule $U$ of the bimodule
$\overline{V^{F}}$ such that $V^{F}=V({\bf C}(F_3))\oplus U$,
$\partial (t)\in V({\bf C}(F_3))$ for each $t\in B_1$ and
$\partial (u)\in V({\bf C}(F_3))$ for each $u\in V({\bf C}(F_3))$.
Hence it follows from the definition of  the representations
of a box that 
${\bf Ver}\,{\rm rep}\,({\cal C}(F_3))={\bf Ver}\,{\rm rep}\,({\cal A}^{F})$.
Then it is ease to see that 
$(GF^{\ast})(M)\cong P(M)^{\bullet}$
for each indecomposable representation $M$ of ${\bf C}(F_3)$,
where $G$ is as in Proposition 2.
\end{proof}

\begin{definition} Consider the bunch of semi-chains ${\bf C}(F_3)$. 

\begin{itemize}
\item Denote by $S$  the set of all  usual strings $w$
and special strings $(w,k)$, where
$w=w_0r_1w_1r_2w_2\cdots$ $r_mw_m$,
such that one of the following conditions hold:

(a) $w_0=q[t-1]$, $r_1= \sim$,
$w_1=z[t]$ for some $t\in \Z$ and if $w_k\in \{
x[i+1], y[i], z[i], r[i-1], p[i], q[i-1]\}$ then $i\geq t$;

(b) $w_m=q[t-1]$, $r_m= \sim$,
$w_{m-1}=z[t]$ for some $t\in \Z$ and if $w_k\in \{
x[i+1], y[i], z[i], r[i-1], p[i], q[i-1]\}$ then $i\geq t$;

(c) $w_0=r[t-1]$, $r_{1}= \sim$,
$w_{1}=x[t+1]$, $w_2\ne r[t+1]$ for some $t\in \Z$ and if $w_k\in \{
  x[i+1], y[i], z[i], r[i-1], p[i], q[i-1]\}$ then $i\geq t$;

(d) $w_m=r[t-1]$, $r_m= \sim$,
$w_{m-1}=x[t+1]$, $w_{m-2}\ne r[t+1]$ for some $t\in \Z$ 
and if $w_k\in \{
x[i+1], y[i], z[i], r[i-1], p[i], q[i-1]\}$ then $i\geq t$.

\item Denote by $\Psi$  the set of all 
$M\in {\bf Ver}\, {\rm rep}\,({\bf C}(F_3))$ which 
correspond to the elements from $S$.

\end{itemize}
\end{definition}

\begin{theorem}
The projective complexes
$P(M)^{\bullet}, M\in {\bf Ver}\, {\rm rep}\,({\bf C}(F_3))$ and 
$\b(P(M)^{\bullet})^{\bullet}$, $M\in \Psi$, 
constitute an exhaustive list of pairwise non-isomorphic
indecomposable objects of ${\bf D}^{b}(F_3)$.
\end{theorem}

\begin{proof}
It is ease to see that 
${\rm Ker}\, p(\a_1)={\rm Ker}\, p(\b_2)=0$,
${\rm Ker}\, p(\a_2)={\rm Ker}\, p(\b_1)={\rm Ker}\, p(\b_1\a_1)$, and
${\rm Ker}\, p(\b_1)$
has the following minimal projective resolution:

$$\cdots\rightarrow {\bf P}_1\oplus {\bf P}_3
\stackrel{\f}{\rightarrow}  {\bf P}_2\stackrel{\psi}{\rightarrow}
{\bf P}_1\oplus {\bf P}_3 
\stackrel{\f}{\rightarrow}  {\bf P}_2\stackrel{\psi}{\rightarrow}
{\bf P}_1\oplus {\bf P}_3 
\stackrel{\f}{\rightarrow}  {\bf P}_2\stackrel{\psi}{\rightarrow}
{\bf P}_1\oplus {\bf P}_3 
\rightarrow
{\rm Ker}\, p(\b_1) 
\rightarrow 0,$$   

\noindent where 
$$\f=
\left( \begin{array}{ll}
p(\a_1)&p(\b_2)
\end{array}
\right)^{T},
\psi=
\left( \begin{array}{ll}
p(\b_1)&p(\a_2)
\end{array}
\right).$$

Straightforward calculation, which we omit, shows that
$${\cal X}(F_3)=\{P(M)^{\bullet}\,|\,M\in \Psi\}.$$
Hence the theorem follows from Corollary 1 and Proposition 3.
\end{proof}

\begin{center} 
{\bf\large Acknowledgments} 
\end{center} 
The part of this research was done during the visit of the first author to 
University of S\~ao Paulo. The final version was prepared during the visit of the second author to the Fields Institute.  
The financial support of FAPESP and the Fields Institute is gratefully acknowledged. 
We are indebted 
to D.Melville for helping us out with Young diagram computations.

\vspace{1cm}
 
\noindent
Viktor Bekkert, 
Universidade Federal do Rio Grande do Norte,
Departamento de Matem\'a\-tica,
CCET, Campus Universit\'ario, Lagoa Nova,
CEP 59072-970, Natal - RN, Brasil,
e-mail: {\tt bekkert\symbol{64}ccet.ufrn.br}
\vspace{0.4cm}

\noindent
Vyacheslav Futorny, 
Universidade de S\~ao Paulo, 
Instituto de Matem\'atica e Estatistica, 
Caixa Postal 66281 - CEP 05315-970, 
S\~ao Paulo, Brasil, 
e-mail: {\tt futorny\symbol{64}ime.usp.br}


\begin{thebibliography}{16}
\bibitem[B]{B} V. V. Bondarenko, 
{\em Representations of bundles of semi-chained sets and their applications},
Algebra i Analiz {\bf 3} (5) (1991), 38--61;
English transl.: St. Petersburg Math. J. {\bf 3} (1992), 973--996.

\bibitem[BM]{BM} V. Bekkert, H. Merklen, 
{\em Indecomposables
in derived categories of gentle algebras}, 
preprint RT-MAT 2000-23, University of S\~ao Paulo (2000), 18p;
to appear in Algebras and Representation Theory.

\bibitem[BGG]{BGG} I. Bernstein, I. Gelfand, S. Gelfand, 
{\em On certain category of ${\cal G}$-modules}, 
Funkt. Anal. i Ego Prilozh. {\bf 10} (1976), 1--8.

\bibitem[CB]{CB} W. W. Crawley-Boevey,
{\em Functorial filtrations II: Clans and the Gelfand problem},
J. London Math. Soc. (2) {\bf 40} (1989), 9--30.

\bibitem[D]{D} Yu. A. Drozd, 
{\em Tame and wild matrix problems}, 
In:  Representations and quadratic forms, Kiev (1979), 39--74;
English transl.: AMS Translations, {\bf 128} (1986), 31--55.

\bibitem[De]{De} B. Deng,
{\em On a problem of Nazarova and Roiter},
Comment. Math. Helv. {\bf 75} (2000), 368--409.

\bibitem[DEMN]{DEMN} S. Doty, K. Erdmann, S. Martin, D. Nakano, 
{\em Representation type of Schur algebras}, 
Mat. Z. {\bf 233} (1)  (1999), 137--182.

\bibitem[DG]{DG} Yu. A. Drozd, G-M. Greuel,
{\em On the classification of vector bundles on projective curves},
preprint (1999).

\bibitem[DN]{DN} S. Doty, D. Nakano, 
{\em Semisimplicity of Schur algebras}, 
Math. Cambridge Phil. Soc. {\bf 124} (1998), 15--20.

\bibitem[DNP1]{DNP1} S. Doty, D. Nakano, K. Peters, 
{\em Infinitesimal Schur algebras}, 
Proc. London Math. Soc. {\bf 72} (3) (1996), 588--612.

\bibitem[DNP2]{DNP2} S. Doty, D. Nakano, K. Peters, 
{\em Infinitesimal Schur algebras 
of finite representation type}, 
Quart. J. Math. Oxford {\bf 48} (1997), no.191, 323--345.

\bibitem[DNP3]{DNP3} S. Doty, D. Nakano, K. Peters, 
{\em Polynomial representations of Frobenius 
kernels of $GL_2$}, 
Contemporary Math. {\bf 194} (1996), 57--67. 

\bibitem[E]{E} K. Erdmann, 
{\em Schur algebras of finite type},
Quart. J. Math. Oxford {\bf 44} (2) (1993), 17--41. 

\bibitem[G]{G} J. Green, 
{\em Polynomial Representations of $GL_n$}, 
Lecture Notes in Math. 
{\bf 830}, Springer-Verlag, New-York (1980).

\bibitem[GK]{GK} Ch. Geiss, H. Krause, 
{\em On the notion of derived tameness}, 
preprint (2000), www.matem.unam.mx/~christof/preprints/derived.ps.

\bibitem[GR]{GR} P. Gabriel, A. Roiter, 
{\em Representations of finite-dimensional algebras}, 
Algebra VIII, Encyclopedia of Math. Sc. {\bf 73}, Springer (1992).

\bibitem[H]{H} D. Happel, 
{\em Triangulated categories in the representation
theory of finite-dimensional algebras}, 
Cambridge University Press, Cambridge (1988).

\bibitem[Har] {Har}  R. Hartshorne, 
{\em Residues and Dualities},
Springer LNM {\bf 20} (1966). 

\bibitem[HR]{HR} D. Happel, C. M. Ringel, 
{\em The derived category of a tubular algebra}, 
Springer LNM {\bf 1177} (1984), 156--180.

\bibitem[JK]{JK} G. James, A. Kerber, 
{\em The representation theory of the symmetric group}, 
Addison-Wesley, 1981.

\bibitem[KZ]{KZ} S. K\"onig, A. Zimmermann, 
{\em Derived equivalences for group rings}, 
Springer LNM {\bf 1685} (1998).

\bibitem[P]{P} B. Parshall, 
{\em Finite dimensional algebra and algebraic groups}, 
Contemporary Mathematics {\bf 82} (1989), 97--114.

\bibitem[Ri]{Ri} C. M. Ringel, 
{\em Tame algebras and integral quadratic forms}, 
Springer LNM {\bf 1099} (1984).

\bibitem[Ro]{Ro} A. V. Roiter, 
{\em Matrix problems and representations of BOCS's}, 
Springer LNM {\bf 831} (1980), 288--324.

\bibitem[U]{U} L. Unger, 
{\em The concealed algebras of minimal wild hereditary algebras}, 
Bayreuther Mathematishe Schriften {\bf 31} (1990), 145--154.

\bibitem[V]{V} D. Vossieck, 
{\em The algebras with discrete derived category},
Preprint (2000).

\bibitem[W]{W} Ch. A. Weibel,
{\em An introduction to homological algebra},
Cambridge Studies in Advanced Mathematics {\bf 38},
Cambridge University Press (1994).

\end{thebibliography}
\end{document}